\documentclass[english,11pt]{article}
\usepackage[T1]{fontenc}
\usepackage{babel}
\usepackage{amsmath,amsfonts,amsthm}
\usepackage{xcolor}
\usepackage{graphicx}
\usepackage{pdfsync}
\usepackage{hyperref}
\hypersetup{linktocpage}
\usepackage[all]{xy}
\hypersetup{colorlinks,citecolor=blue,filecolor=black,linkcolor=blue,urlcolor=black}

\usepackage{fancyhdr}

\DeclareMathSymbol{\leqslant}{\mathalpha}{AMSa}{"36} 
\DeclareMathSymbol{\geqslant}{\mathalpha}{AMSa}{"3E} 

\newcommand{\Star}{\times\kern-11pt+}

\newcommand{\Abar}{{\backslash\kern-8pt A}}

\begin{document}
\renewcommand{\refname}{References in chronological order}

\begin{center}{\Large\sc
Some Historical Aspects of \vspace{.5cm}

Error Calculus by Dirichlet Forms}\\ \vspace{.5cm}

Nicolas Bouleau
\end{center}

\noindent{\bf Abstract.}

We discuss the main stages of development of the error calculation since the beginning of XIX-th century by insisting on what prefigures the use of Dirichlet forms and emphasizing the mathematical properties that make the use of Dirichlet forms more relevant and efficient. The purpose of the paper is mainly to clarify the concepts. We also indicate some possible future research.\\

\noindent{\bf I. Introduction.}  There are several kinds of error calculations which have not followed the same historical development. The error calculus by Dirichlet forms that we will explain and trace the origins has to be distinguished from the following calculations:

a) The calculus of roundoff errors in numerical computations which appeared far before the representation of numbers  in floating point be implemented on computers, and which possesses its specific difficulties. It has been much studied during the development of the numerical analysis for matrix discretization methods  (cf. Hotelling \cite{hotelling}, Von Neumann  \cite{vonneumann}, Turing \cite{turing}, Wilkinson \cite{wilkinson}, etc.);  

b) The global evaluation of deterministic errors such as the  interval calculus (cf. Moore \cite{moore}, etc..). Some works of Laplace are related to this approach and also the paper of Cauchy  \cite{cauchy}; 

c) The calculus of finite probabilistic errors where the errors are represented by random variables, which has been used by a very large number of authors to begin an argument and then, often, modified by supposing the errors to be small or gaussian in order to be able to pursue the calculation further (cf. Bienaym\'e \cite{bienayme}, Birge \cite{birge}, Bertrand \cite{bertrand}, etc.) because the computation of image probability distributions is concretely inextricable what, in the second half of the XX-th century, justified the development of simulation methods (Monte-Carlo and quasi-Monte-Carlo).\\ 

The error calculus by Dirichlet forms assumes the errors to be both small, actually  infinitesimal, and probabilistic. These two characteristics imply a peculiar differential calculus for the propagation of errors through models. As we will see the part of the calculation related to what is called today {\it the squared field operator} or more often {\it the carr\'e du champ operator}, is   ancient and dates back to the turn of the XVIII-th and XIX-th centuries in connection with the birth of the least squares method. Let us note, however, that our purpose is not to make a history of the method of least squares, broad topic that would lead to decline all the benefits of optimization in $L^2$ and its developments in statistics and analysis. I refer in this regard to the historical work of Kolmogorov and Yushkevich  \cite{kolmogorov-yushkevich},  also to the book of  Pearson \cite{pearson}, and to the article of  Sheynin \cite{sheynin} not always clear from a mathematical point of view probably because of an intrinsic ambiguity of the thought of the authors of the turn of the XVII-th and XIX-th centuries, but extremely well documented.  Unfortunately it does not address at all  the propagation of errors. \\

Let us note also in this introduction the very important phenomenon of  {\it dichotomy of small errors} which allows to clarify the validity domains of the probabilistic and deterministic approaches. When we are concern with {\it small} errors that means mathematically that we are in an approximation procedure where -- in principle -- we are able to make the errors vanish. In such a situation the description of the asymptotic mechanism of the propagation of infinitely small errors is different according to the respective size of the variance and the bias of the error.

Three cases appear in the limited expansions.  1) If the variance is negligible with respect to the bias, then this property will persist by deterministic computations and it is enough for the asymptotic calculus to perform a first order differential calculus, i.e. a classical sensitivity calculation.  2) If the variance is of the same order of magnitude as the bias, then the calculus has to be a first order differential calculus for the variances which does not involve the bias, and the calculus for the bias is a second order differential calculus involving both biases and variances.  3) If the bias is negligible with respect to the variance then from the first calculation we return to the case 2). In these two last cases the propagation formulae are the following for a scalar erroneous quantity   $X$:
\begin{equation}\hspace{-.5cm}\label{propagation}
\left.
\begin{array}{rl}
\mbox{bias of error on }f(X)&=(\mbox{bias of error on }X) f^\prime(X)+\frac{1}{2}(\mbox{var of error on }X)f^{\prime\prime}(X)\\
\mbox{var of error on }f(X)&=(\mbox{var of error on }X)f^{\prime 2}(X)
\end{array}
\right\}
\end{equation}
We refer to  \cite{bouleauRIMS} for more details and typical examples. Situations like case 1) are called {\it weakly stochastic}, and situations like 2) or 3) are called {\it strongly stochastic}. Let us mention that often in practice, in engineering for instance,  we are not able to  control the nature of the errors.  Errors on the data in modelling are {\it exogenous}, we know few where they come from. It is wise to think according to the case of strongly stochastic errors, especially to take in account the randomness of the errors through the non-linearities of the model.

This is important because, by extending the ideas of Poincar\'e on the errors  \cite{poincare}, it is possible to see that measurements done with graduated instruments are always  {\it strongly stochastic} (because of the error on the choice of the nearest graduation cf. \cite{bouleauhal} \cite{bouleau3}, \cite{bouleau2006b}). This is related to the theory of {\it arbitrary functions} (cf. von Kries \cite{vonkries}, Fr\'echet \cite{frechet}, Borel \cite{borel}, Hostinsk\'y \cite{hostinsky1}, Hopf \cite{hopf1}, \cite{hopf2}, \cite{hopf3}).\\

In summary we can distinguish
\begin{table}[h] 
\begin{center}
\noindent\begin{tabular}{p{2cm}||p{2.5cm}|p{2.5cm}|p{2cm}|}
&\multicolumn{2}{c|}{}& \\ 
&\multicolumn{2}{c|}{Infinitesimal errors}& \multicolumn{1}{c|}{Finite errors}\\ 
&\multicolumn{2}{c|}{}&\\ 
 \hline\hline
  & \multicolumn{2}{c|}{Deterministic sensitivity analysis :}& \\
\multicolumn{1}{c||}{Deterministic} &\multicolumn{2}{c|}{derivation with respect to}& \multicolumn{1}{c|}{Interval calculus,}\\
\multicolumn{1}{c||}{approaches} & \multicolumn{2}{c|}{the parameters of the model}& \multicolumn{1}{c|}{sup norm}\\
&\multicolumn{2}{c|}{}&\\ \hline
&\multicolumn{2}{c|}{Error calculus using}&\\
\multicolumn{1}{c||}{Probabilistic} & \multicolumn{2}{c|}{Dirichlet forms}& \multicolumn{1}{c|}{Probability theory :}\\
\multicolumn{1}{c||}{approaches}&\hspace{.5cm}first order&\hspace{.3cm}second order & \multicolumn{1}{c|}{images of}\\ 
\multicolumn{1}{c||}{}& calculus only &calculus with
& \multicolumn{1}{c|}{probability}\\ 
\multicolumn{1}{c||}{}&dealing with&
variance\index{variance}s and & \multicolumn{1}{c|}{distributions}\\
\multicolumn{1}{c||}{}&variances& bias\index{bias}es& \multicolumn{1}{c|}{}\\\hline
\end{tabular}
\caption{Main categories of error calculi}
\end{center}
\end{table}
\newpage
\noindent{\bf II. Gauss inventor of the  carr\'e du champ operator ?} 

\noindent Referring to the turn of the XVIII-th and the XIX-th centuries   we see that two wakes  were clearly drawn among researchers in matter of error\;: that of Laplace and that of Gauss.

Let us begin by Gauss whose works interess us particularly here. The  current he initiates gives a fundamental place to his {\it law of errors}\;: 
 assuming that, after several independent measures $x_i$ the arithmetic
 average  $\frac{1}{n}\sum_{i=1}^n x_i$ is the best value to take in account, he showed, with some
 additional assumptions, 
that the errors follow necessarly a normal law and the arithmetic average is both the most likely value and the one given
by the least squares method.

Gauss tackled this question in the following way. First he admits -- and this idea will be kept in the error theory with Dirichlet forms
 -- that the quantity to be measured is random. It can vary in the scope of
 the measurement device following an {\it a priori} law. In modern language, let $X$ be this random variable, $\mu$ its law.
The results of the measurement operations are other random variables $X_1,\ldots, X_n$ and Gauss assumes

a) the conditional law of $X_i$ given $X$ to be of the form $\mathbb{P}\{X_i\in E|X=x\}=\int_E \varphi(x_1-x)\;dx_1$,

b) the variables $X_1,\ldots,X_n$ to be conditionally independent given $X$.

Then he is easily able to compute the conditional law of $X$ given the results of measure,
it has a density with respect to $\mu$ and writing this density is maximal at the arithmetic average,
 he gets the relation
$$\frac{\varphi^\prime(t-x)}{\varphi(t-x)}=a(t-x)+b$$
hence 
$$\varphi(t-x)=\frac{1}{\sqrt{2\pi\sigma^2}}
\exp(-\frac{(t-x)^2}{2\sigma^2}).$$
\indent It is probably in the course {\it Calcul des Probabilit\'{e}s} of Poincar\'{e} \cite{poincare} at the end of the XIXth century
 that Gauss' argument is the most clearly explained because Poincar\'{e} tries both to explicit 
all hypotheses and to generalize the proof\footnote{It is about this `law of errors' that Poincar\'{e}
writes ``everybody believes in it because experimenters imagine it is a theorem of mathematics
and mathematicians it is an experimental fact".}. He studies the case where the conditional law of
 $X_i$ given $X$ is no more $\varphi(y-x)\,dy$ but of the more general form $\varphi(y,x)\,dy$. That gives him 
a way to explain the seeming paradox of {\it the error permanency}: "with a meter divided in millimeters,
he writes, as often the measures be repeated, never will a length be determined up to a millionth 
of millimeter". This phenomenon is well known by physicists who noted of course that during 
the whole history of experimental sciences never a quantity has been precisely measured with rough instruments, cf \cite{klein}.
The discussion leads Poincar\'{e} to suggest that measurements could be independent whereas the errors be not,
 when done by the same
instrument. He doesn't develop a new mathematical formalism for this idea but emphasizes on the advantage of assuming 
small errors because then the argument of Gauss giving the normal law
becomes compatible with non linear changes of variables and can be performed through dif\-ferential calculus.
 It is the question of the
{\it error propagation}.\\

It is in his  {\it Theoria Combinationis Observatonum Erroribus Minimis Obnoxiae} published in  1823 that Gauss details his ideas about the errors propagation. 

In the introduction, he cites Laplace and discusses the merits of reasonning with repeated observations or with observations immediately erroneous. Behind this discussion is the fact that Laplace gave the first analytical proof of the central limit theorem, and that Gauss  intends to assert the interests of his demonstration that if the arithmetic mean is taken as the correct value then the law is necessarily normal, that he replaces in a more general approach of some kind of error calculus in an extended meaning.  This dicussion  is deepened in his section  17.

Gauss states the main problem in the following way:
{\it Given a quantity $U = F(V_1,V_2,V_3,\ldots)$ function of the erroneous quantities $V_1, V_2, V_3,\ldots$, compute the potential quadratic error to expect on $U$ with the quadratic errors $\sigma_1^2, \sigma_2^2, \sigma_3^2, \ldots$ on $V_1, V_2, V_3, \ldots$ being known and assumed small and independent.}

His answer is the following formula :
\begin{equation}\label{propagationgauss}
\sigma_U^2=(\frac{\partial F}{\partial V_1})^2\sigma_1^2+(\frac{\partial F}{\partial V_2})^2\sigma_3^2+(\frac{\partial F}{\partial V_3})^2\sigma_3^2+\cdots
\end{equation}
He also provides the covariance between an error on $U$ and an error on another function of the $V_i$'s.
Formula (\ref{propagationgauss}) displays a property which makes it much to be preferred  to other formulae encountered in textbooks throughout the XIXth and XXth centuries. It features a coherence property. With a formula such as
\begin{equation}\label{ugly}\sigma_U=|\frac{\partial F}{\partial V_1}|\sigma_1+|\frac{\partial F}{\partial V_2}|\sigma_3+|\frac{\partial F}{\partial V_1}|\sigma_3+\cdots
\end{equation}
errors may depend on the way in which the function $F$ is written. These "ugly" formulae remain for instance in \cite{baret}, \cite{vessiot}.
Today we can understand that this difficulty does not arise with Gauss' calculus thanks to its connection with the theory of Dirichlet forms. Introducing the differential operator 
\begin{equation}L=\frac{1}{2}\sigma_1^2\frac{\partial^2}{\partial V_1^2}+\frac{1}{2}\sigma_2^2\frac{\partial^2}{\partial V_2^2}+\frac{1}{2}\sigma_3^2\frac{\partial^2}{\partial V_3^2}\cdots
\end{equation}
and supposing the functions to be smooth, we remark that formula (\ref{propagationgauss}) can be written as 
\begin{equation}\sigma_U^2=L(F^2)-2FL(F)\end{equation}
and coherence follows from the transport of a differential operator by an application. Today these intrinsic properties are understood also by the link with the stochastic calculus on manifolds  and particularly the differential calculus of second order \cite{malliavin}, \cite{meyer81}, \cite{bouleau2006a}.

The errors on $V_1, V_2, V_3, \ldots$ are not necessarily supposed to be independent nor constant and may depend on $V_1, V_2, V_3,\ldots$, that gives formulas obtained by, what we call today, polarisation.

For the proof of  (\ref{propagationgauss}) Gauss performs -- as will do most of the applied textbooks thereafter -- a computation supposing errors are {\it infinitely small}  quantities and may replace differentials  $dV_1, dV_2$, etc. But even Gaussian random variables have non compact support. There exists therefore a small probability that the errors be large. This is a difficulty that can be only  treated with more precise notions of convergence or with the theory of Dirichlet forms, the proof of Gauss will then appear as a computation of what we call a {\it gradient in the sense of Dirichlet forms} often denoted by the term {\it sharp} (cf \cite{bouleau3}), we will return to this point for clarification below in part V.

If we read the sections  19 and  20 of his treatise, we may ask whether the question posed by Gauss is not -- in germ -- the idea of a quadratic form which would be the ancestor of a Dirichlet form. 

In some sense Gauss is on a product error structure  (he assumes the errors of observation to be independent -- in the sake of simplicity he says, because he is able to treat the general case). He obtains that the error on a quantity is given by a random quadratic function that he gets by local linearization with his change of variable formula, and whose expectation is the mean quadratic error : this random function is the  "carr\'e du champ". 
I do not say that Gauss had the idea of the carr\'e du champ, but I say that his direction of research has something to do with Dirichlet forms what is not philosophically surprising since the landscape of potential theory was for him familiar.\\

Let us say now some words of the wake left by Laplace. By studying (cf \cite{laplace}, \cite{laplace2}, \cite{laplace3} and Cauchy \cite{cauchy}) the mean of the positive errors and that of negative errors then taking the mean anew,  Laplace, actually studies the first absolute moment of the error 
$$\mathbb{E}[\,|\mbox{error}|\,]$$
instead of its variance. It is what seems for him the most natural. He has among available tools to make calculation the method of characteristic functions that he perfectly masters, and for this he must choose hypotheses on the law of this random variable $|\mbox{error}|$. But for the propagation he faces the usual difficulty of intractable calculations. He is then led to assume small and Gaussian errors following the argument provided by the central limit theorem when observations are repeated.

But when the errors are Gaussian the moments of order 1 and 2 are linked and we have, if $V$ is $N(0,\sigma)$ distributed, 
$$\mathbb{E}[\,|V|\,]= \sqrt{\frac{2}{\pi}}\;\sigma$$
so that, if the errors are small, the relation between variances is {\it equivalent} to the following relation for the first absolute moment : 
$$\mathbb{E}[\,|e_F|\,]=\sqrt{F_1^{\prime 2}(\mathbb{E}|e_1|)^2+\cdots+F_n^{\prime 2}(\mathbb{E}|e_n|)^2}.$$
We see that for first order moment a formula is obtained very similar to the one of Gauss  (\ref{propagationgauss}) (instead of  "ugly" formulae like (\ref{ugly}) ).\\

In the wake of the works of Gauss on errors we could cite almost all treatises on probability of the XIX-th century, let us quote particularly  Faa di Bruno \cite{faadibruno} who follows narrowly his approach and cites him explicitely  (p38 et 53), also Liagre \cite{liagre}. 

In the wake of  Laplace we can cite Bravi \cite{bravi}, Bravais \cite{bravais} who extends the arguments of Laplace to the multivariate case and introduces the error ellipsoids  (cf \cite{bouleau3} Chap.I), A. Meyer (1857) who in a treatise of very good mathematical level  \cite{ameyer}  opts for the approach of Laplace and argues this position, also the physicist  Airy (1861) \cite{airy} who takes the whole theory of probability starting only from the treatise of Laplace.

Belonging clearly to both wakes are  Hagen \cite{hagen}, Bienaym\'e \cite{bienayme} and also  Biver \cite{biver} who, in 1853 introduces a  general {\it cost functional} to be minimized to manage errors.\\

The propagation formula is encountered also in statistics if one is interested in  the propagation of the Fisher information  \cite{fisher}, \cite{fisher2}. If $\mathbb{P}_\theta$, $\theta\in \mathbb{R}^d$, is an indexed family of probability measures satisfying the conditions of regular models, having a density  $f(.,\theta)$  with respect a measure  $\mathbb{Q}$, the Fisher information related to the parameter   $\theta$ is the matrix
$$J= 4\int (\frac{d\sqrt{f(.,\theta)}}{d\theta})(\frac{d\sqrt{f(.,\theta)}}{d\theta})^t\,d\mathbb{Q}.$$
If instead of the parameter  $\theta$ we consider the parameter  $\varphi=g(\theta)$ where $g$ is smooth and injective, we find that the inverse matrix  $\Gamma(\theta)= J^{-1}(\theta)$ that represents  {\it a precision} on $\theta$, is transported following Gauss formula. 
$$\Gamma_{ij}(\varphi)=\sum_{k\ell}\frac{\partial g_{ij}}{\partial \theta_k}\frac{\partial g_{ij}}{\partial \theta_\ell}\Gamma_{k\ell}(\theta)$$
see on this topic   \cite{bouleauchorro} et Chorro \cite{chorro}. This link with statistics gives a new interpretation of the carr\'e du champ operator in terms of accuracy of a statistical estimate, there are no small errors more, only information and its inverse the precision. In this interpretation the Fr\'echet-Cramer-Rao inequality says that efficient estimates give the best precision. 

We understand then that obtaining bounds for the propagation of the Fisher information is a means to express the regularity of a carr\'e du champ matrix  cf. Villani \cite{villani} p824 {\it et seq.}\\

\noindent{\bf III. Why should we ask the quadratic form to be closed?}
\\*The Dirichlet form strictly speaking appeared in potential theory in the classical case of the Laplacian operator far before the thought structures of the Hilbert spaces be available, it was simply a quantity whose minimum value was looked for. Dirichlet uses it \cite{dirichlet} as soon as 1846, Thomson \cite{thomson} and Riemann \cite{riemann} near after, we refer to the historical sudies on the potential theory itself of  Brelot \cite{brelot},   Temple \cite{temple}, Ancona \cite{ancona} among others on these developments.

After the works of Henri Cartan \cite{cartan} and the thesis of  Jacques Deny \cite{deny}, \cite{deny-bis} it is during the collaboration between  Arne Beurling and Jacques Deny \cite{beurling-deny1}, \cite{beurling-deny2}, \cite{deny1}
that Hilbertian methods have taken a new importance especially by this discovery that the fact that contractions operate on a closed positive quadratic form is a necessary and sufficient condition for the associated contraction semi-group  be positive on positive functions, i.e.$\!$ be the transition probability semi-group of a Markov process. 

The study of these questions under the angle of symmetric Markov processes developed quickly during the 1970s, Silvestein \cite{silverstein}, Le Jan \cite{lejan}, (see also the references of the contributions to the  Colloque en l'honneur de J. Deny, Orsay 1981), what allows Masatoshi Fukushima to discover these famous continuous additive functionals with zero quadratic variation but not finite variation \cite{fukushima0}.

The book of  Fukushima \cite{fukushima}, has been a reference for researchers and the starting point of several works. The notion of Dirichlet form appeared there so clearly -- a primary notion somehow -- that it came immediately to the mind that this notion could accompany any probability space from where the links with Malliavin calculus and error calculus.\\

From the point of view of error theory, the idea of imposing a priori that the quadratic form (the expectation of the carr\'e du champ)  be closed is a major step.  The error calculus of Gauss contains the limitation of supposing that both the function $F$ and the random variables $V_1, V_2, V_3, \ldots$ are analytically known. In probabilistic modelling however, we are often confronted to situations in which all the random variables, functions and covariances matrices are {\it given by limits}. For such situations, a means of extension thereby becomes essential.

Let us suppose that  the quantities be defined on the probability space $(\Omega, \mathcal{A},\mathbb{P})$. The extension tool lies in the following : we assume that if $X_n\rightarrow X$ in $L^2(\Omega, \mathcal{A},\mathbb{P})$
and if the variance of the error $\Gamma(X_n-X_n)$ on $X_n-X_n$ can be made as small as we
wish in $L^1(\Omega, \mathcal{A},\mathbb{P})$ for $m, n$ large enough, then the variance of the error $\Gamma(X_n-X)$ on $X_n-X$ goes to zero in $L^1$.
This idea can be interpreted as a reinforced coherence principle, it means that the error on $X$ is attached to $X$ and furthermore, if the sequence of pairs ($X_n$, error on $X_n$) converges suitably, it converges necessarily to a pair ($X$, error on $X$).

The main benefit of the extension tool is that error theory based on Dirichlet forms extends to the infinite dimension, which allows  error calculus on stochastic processes (Brownian motion, Poisson random measures, diffusions defined by stochastic differential equations). But also the calculus itself becomes more flexible : it allows now Lipschitzian changes of variables.\\

From a historical perspective, we can say that the main characteristic of the XX-th century with respect to XIX-th century, already so extraordinarily developed in mathematics, is the use of functional spaces (where points represent functions) and especially {\it complete} functional spaces like {\it Hilbert spaces, Banach spaces},  that permit to handle objects defined by limits. This idea has been the origin of gigantic progress in all domain of analysis. The closedness of Dirichlet forms allows to install an error calculus in most situations of classical analysis and stochastic  analysis.

We arrive then to a kind of enriched error calculus, where any development -- in particular asymptotic theorems -- may be accompanied by errors. For instance the famous Donsker theorem on the approximation of the Brownian motion by a random walk extends in term of erroneous random walk which yields naturally the Ornstein-Uhlenbeck structure on the Wiener space cf \cite{bouleauD}.\\

\noindent{\bf IV. Dirichlet form generated by an approximation}
\\*If we want to perform very precise error calculations -- this is not an oxymoron\;! -- we must be concern by the biases.  Among the engineers who ask for a great precision from the observations, there are the geodesists and it seems that it was in this domain that has been used for the first time the fact that the random nature of errors imply a bias.  No doubt everyone has long known that the image of the expectation  for a non-linear mapping is not the expectation  of the image,  but geodesists write down differential error calculi with biases. This is explicit  in Elkins (1950) \cite{elkins1}, \cite{elkins2}, \cite{elkins3}, Jeudy (1988) \cite{jeudy}, Hangos-Leisztner (1989) \cite{hangos}, Teunissen \cite{teunissen}\cite{teunissen2}, Coleman-Steele \cite{coleman}.

In Section  17 of the {\it Theoria Combinationis} Gauss writes "It is necessary to warn here that in the following research, it will issue only random errors reduced of their constant part. It is the observer's responsibility to remove the causes of constant error carefully. We reserve for another opportunity to examine the case where the observations are affected to an unknown error, and we will discuss this issue in another memory" maybe he would have encountered the question of the bias\ldots

As we will see in the framework of Dirichlet forms the bias is represented by the generator of the semi-group associated with the form. The exact formula for the propagation of the bias is the following   (see \cite{bouleau-hirsch2} p42 Exercise 6.2 and the precise hypotheses required for this formula)
$$A(F(f))=\sum_{i=1}^n F_i^\prime(f)Af_i+\frac{1}{2}\sum_{i,j}F_{i,j}^{\prime\prime}\Gamma(f_i,f_j).$$
A similar expression is used by the geodesists.\\

But the bias is a delicate notion because of some latitude in its definition. Let us resume what says the theory. If we consider two random variables $Y$ and $Y_n$ close together, the asymptotic behaviours of 
 $$\mathbb{E}[(\phi(Y_n)-\phi(Y))\chi(Y)]$$
 and of 
 $$\mathbb{E}[(\phi(Y_n)-\phi(Y))\chi(Y_n)]$$
 where $\phi$ and $\chi$ are test functions, are generally different. As a consequence several bias operators have to be distinguished (cf. \cite{bouleau2006}) :
 
The asymptotics of $\mathbb{E}[(\varphi(Y_n)-\varphi(Y))\chi(Y)]$ yields the theoretical bias operator $\overline{A}$,

\noindent the asymptotics of $\mathbb{E}[(\varphi(Y)-\varphi(Y_n))\chi(Y_n)]$ yields the practical bias operator noted $\underline{A}$

\noindent that of $\mathbb{E}[(\varphi(Y_n)-\varphi(Y))(\chi(Y_n)-\chi(Y))]$ gives the symmetric bias operator $\widetilde{A}$,

\noindent and eventually $\mathbb{E}[(\varphi(Y_n)-\varphi(Y))(\chi(Y_n)+\chi(Y))]$ provides the singular bias operator $\Abar$.

These operators are related thanks to the relations
$$\widetilde{A}=\frac{\overline{A}+\underline{A}}{2}\quad\quad\Abar=\frac{\overline{A}-\underline{A}}{2}.$$
The symmetric bias operator satisfies
$$<\widetilde{A}[\varphi],\chi>_{L^2(\mathbb{P}_Y)}=<\varphi,\widetilde{A}[\chi]>_{L^2(\mathbb{P}_Y)}$$
under natural hypotheses it is indeed the generator of a Dirichlet form ${\cal E}$ which possesses a carr\'e du champ operator $\Gamma$. It may be shown that this Dirichlet form is local iff the asymptotics (with the same weight) of 
$\mathbb{E}[(\varphi(Y_n)-\varphi(Y))^4]$ vanishes.

If the asymptotics of both variances $\mathbb{E}[(\varphi(Y_n)-\varphi(Y))^2\psi(Y)]$ and $\mathbb{E}[(\varphi(Y_n)-\varphi(Y))^2\psi(Y_n)]$ coincide -- what is usually the case -- the singular bias operator is a first order operator in the sense that it satisfies on the test functions
$$\Abar[\varphi\chi]=\Abar[\varphi]\chi+\varphi \Abar[\chi].$$

Let us take an example. A very frequent  situation in probability theory involve  a triplet of real random variables $(Y,Z,T)$, a real random variable $G$ independent of $(Y,Z,T)$ centered with variance one, and
the approximation $Y_\varepsilon$ of $Y$ given by
\begin{equation}
Y_\varepsilon=Y+\varepsilon Z+\sqrt{\varepsilon}TG.
\end{equation}
In that case the operator $\overline{A}$ may be shown to be given by

$$\overline{A}[\varphi](y)=\mathbb{E}[Z|Y\!=\!y]\;\varphi^\prime(y)+\frac{1}{2}\mathbb{E}[T^2|Y\!=\!y]\;\varphi^{\prime\prime}(y)$$
and the Dirichlet form is given by 
$$\mathcal{E}[\varphi, \chi]=\mathbb{E}[T^2\varphi^\prime(Y)\chi^\prime(Y)]$$
The operator $\widetilde{A}$ depends only on $T$, not on $Z$, $\underline{A}$ is obtained by difference.

A  {\it weakly stochastic} approximation (see Part I above) may now be defined more precisely by the condition $\widetilde{A}=0$.\\

\noindent{\bf V. "Small errors" what does it mean ?}
\\*We can now clarify the question of small errors that has been historically a kind of conceptual blockage.  We have seen that the theory of local Dirichlet forms with a carr\'e du champ operator should be seen as a more accurate and rigorous form of the error calculation developped  in the wake of Laplace and Gauss, allowing an extension of the sensitivity calculation to broader and more difficult situations involving stochastic calculus and Brownian motion.

In the setting of Dirichlet forms we know that the carr\'e du champ operator represents the variance of the error and that the bias of the error can be represented by the symmetric generator (depending on which bias we speak as seen above). But the error itself, what is it\;? The error theory by Dirichlet forms considers that the error is  $X_t-X_0$ where $X_t$ is the symmetric Markov process associated to the Dirichlet form and which is taken as a {\it tangent process to the studied approximation procedure}. The error is not the "sharp". The "sharp" is mathematically a  {\it gradient}. It is a linear operator which restores the carr\'e du champ by taking the square. Actually it is a randomized gradient, for example in the classical case on $\mathbb{R}^2$ it writes $f^\sharp=f_1^\prime \xi_1+f_2^\prime \xi_2$ where $\xi_1$ and $\xi_2$ are auxiliary orthogonal reduced random variables. This looks like an error by the fact that it is random but it is a  tangent vector (a first order differential operator) which acts proportionally to itself.

On the contrary  $X_t$ evolves with a transition probability kernel. The crucial point is that for small  $t$, $X_t-X_0$ is always a  {\it sum of infinitely many quantities with independent sources of randomness}, it is not a quantity which decreases homothetically to itself,  it does not fall under the Taylor formula
$$ f(x+h) = (e^{hD} f)(x)=f(x)+hDf(x)+\cdots+\frac{h^n}{n!}D^nf(x)+\cdots$$
what is relevant for it is the {\it theory of semi-groups of operators}.\\

The errors may receive a rigorous mathematical treatment only if they are thought inside an approximation procedure.  We can say that   {\it throughout the nineteenth century} there was confusion between an error of the type   $hY$ with $h$ a scalar number tending to zero and  an error with accumulation of independence as a Markov process in small time.  In this respect Louis Bachelier \cite{bachelier} with his players of infinitesimal games appears well marking the turning of the two centuries. He will be cited by Kolmogorov \cite{kolmogorov0}.

We find, more in Laplace's writings than in  Gauss' ones (he keeps his distance) the concern that errors are sums of many independent terms. Chapters III and V of the  {\it Th\'eorie Analytique} insistently express this idea. But beyond the central limit theorem, Laplace is obviously not in  position to make a theory where the Taylor formula be replaced by exponentials of operators, i.e. the theory of operators semi-groups. 

In Section 17 of the  {\it Theoria combinationis} Gauss writes
"Laplace, considering the question under a different point of view, showed that this principle [of supposing that the errors follow a normal law] is better than any other choice, for any probability distribution of the errors of observation as soon as the number of observations is very large. But when this number is restricted, the question remains untouched; so that if our hypothetical law [the normal law] is rejected, the least squares method rests better than the others  [e.g. than the use of absolute first moment], by the simple reason that it leads to simpler calculations."

	This remark is relevant, but it shows also  that Gauss somehow underestimates the mathematical importance of the accumulation of small errors. 

Today we know that the infinitesimal reduction of random quantities is well described by the infinitesimal generators of Markov semi-groups which are second order differential operators  (Kolmogorov \cite{kolmogorov0}, Ventsel \cite{ventsel}, Hunt's complete maximum principle  \cite{hunt2} Section 15,  Meyer-Dellacherie \cite{dellacherie4} Theorem XIII 22-24) and have a fractional part  (integro-differential operator) only in the presence of jumps  (case of non local Dirichlet forms, cf below VI \S e). Similarly the stochastic diffential geometry uses second order tangent vectors  (cf \cite{bouleau2006a}, Meyer \cite{meyer81}).

It emerges from this discussion that the "small errors" in Gauss' calculation have to be read as a computation of gradient, what reinforces the importance of the propagation formula under its analytical form  (\ref{propagationgauss}) which is that of the carr\'e du champ operator. \\

\noindent{\bf VI. Trails of research}
\\*a) Obtaining numerical results is often difficult in probability theory because the spaces are high dimensional, often infinite dimensional. Therefore the simulation methods, also called Monte Carlo methods are very usefull thanks also  to their generality,  cf for instance \cite{bouleau-lepingle}. Intuitively the computation of the value of the carr\'e du champ operator  on a random variable in a parametrized modell, say   $X(\omega, \lambda)$, when $\lambda$ is erroneous and the modell is not, is very simple. It is enough to take a {\it cluster} of points in the neighborhood of the value of   $\lambda$ centered on  $\lambda_0$ and with quadratic dispersion  $\sigma^2$ around $\lambda_0$ and then to collect the dispersion of the corresponding values of  $X$,  with $\omega$ being fixed. The dispersion matrix will give the matrix  $\Gamma(X)$ and the discrepancy between  $X(\omega, \lambda_0)$ and the mean of the cluster of the values of  $X$ will yield the bias. 

	If the model is itself erroneous, the method has to be extended with a cluster around  $\omega_0$. This {\it clusters method} for the error calculus has already been used  (Scotti \cite{scotti}) but has not been theoretically studied up to now and many questions arise : optimal number of points in the cluster with respect to its concentration toward  $\lambda_0$, obtention of the mean error  (which is the square root of the Dirichlet form taken on  $X$), use of acceleration by  quasi-Monte-Carlo, etc.
	
	At present theoretical results are in the opposite way : the methods of calculation of densities of probability distributions  (cf for example \cite{parzen}, \cite{bartlett}, \cite{murthy}, \cite{deheuvels}, \cite{elkins3}, \cite{bally-talay2},\cite {silverman}, etc.) may be accelerated  if an error structure is available  (or a Malliavin calculus) see Caballero-Fernandez-Nualart \cite{caballero-fernandez-nualart}, Kohatsu-Higa-Pettersson \cite{kohatsu-pettersson},
Bouleau \cite{bouleauMC}.\\

\noindent b) A purely mathematical enigma remains at present. It is the conjecture of the so-called  {\it Energy Image Density property} or EID. This property is true for any local Dirichlet structure with carr\'e du champ for real valued fonctions (cf Bouleau \cite{bouleau1} Theorem 5 and Corollary 6), (see also the very close method of Davydov-Lifshitz \cite{davydov-lifshits}, \cite{lifshitz}). It has been conjectured in  1986 (Bouleau-Hirsch \cite{bouleau-hirsch1} p251) that (EID) be true for any local Dirichlet structure with carr\'e du champ. This has been proved for the Wiener space equipped with the Ornstein-Uhlenbeck form and for some other Dirichlet structures by   Bouleau-Hirsch (cf \cite{bouleau-hirsch2} Chap.\! II \S 5 and Chap.\! V example 2.2.4), it has been proved by Coquio \cite{coquio} for random Poisson measures on $\mathbb{R}^d$ and  by Bouleau-Denis as soon as the "bottom" space satisfies itself EID. But this conjecture is at present neither proved nor refuted in all generality, it has to be established in each particular framework. (EID) on the Wiener space is now a very frequently used tool to prove existence of density.\\

\noindent c) The error calculus by Dirichlet forms and the mathematically rigourous framework for the carr\'e du champ operator may be used not only for the computation of measurement errors propagation, but to study the effect of   {\it fluctuations} on physical systems. Very often the physicists handle fluctuations as small errors denoted  $\Delta X$ and conduct calculations in the spirit of Gauss calculus. The program is then to write down anew the fluctuation theory and the deviation that it yields for (non linear) measurements thanks to Dirichlet forms. Theoretical advances have been already obtained in the direction cf for example Albeverio-Grothaus-Kondratiev-R\"ockner \cite{a-grothaus-k-r}, Dembo-Deutschel \cite{dembo-deuschel}. 

In order to be concrete, let us look how  L. landau and E. Lifchitz in their famous textbook  \cite{landau-lifchitz} compute the deviation of a stretched string due to the thermal fluctuations  (Chapter XII exercise 8) : 

{\footnotesize Let $\ell$ be the length of the string, $F$ the tension force. Let be a point at the distance $x$ from an end of the string, $y$ its transverse deplacement. To determine $\overline{y^2}$ we have to determine the equilibrium shape of the string for a given deplacement $y$ of the point $x$; there are two segments of straight line betwen the ends and the point $(x,y)$. The work spent for such a deformation of the string is equal to
$$R_{min}=F.(\sqrt{x^2+y^2)}-x)+F.[\sqrt{x^2+y^2}-(\ell-x)] \cong \frac{Fy^2}{2}(\frac{1}{x}+\frac{1}{\ell-x})$$
The mean square is therefore $\quad\overline{y^2}=\frac{T}{F\ell}x(\ell-x).$}

\noindent These results are in perfect accordance with the approach that consists in taking the Ornstein-Uhlenbeck form on the Brownian bridge, because if 
$$X_t=B_t-tB_1=\int_0^1(1_{[0,t]}-t)dB_s$$
then
$$\begin{array}{rl}\Gamma[X_t]&= \sum_0^1(1_{[0,t]}-t)^2ds=t-t^2=t(1-t)\\
\Gamma[X_s,X_t]&= \sum_0^1(1_{[0,s]}-s)(1_{[0,t]}-t)du=s\wedge t-st= s(1-t)\end{array}$$ for $s<t$.\\

\noindent d) If we have an error structure, we have also a capacity theory associated with the Dirichlet form. This allows a refinement of the almost everywhere computations (see Fukushima \cite{fuku-deny} and works in the wake). This may be related with the error theory under the following regard\,: this  helps to understand that during an approximation procedure some things {\it are not seen}.

The fact that there are errors, hence randomness, in the approximation procedure  {\it erases} some features of reality.  In our framework, the procedure is replaced by a Markov process $X_t$ which is in a sense osculating when $t\rightarrow 0$. But the paths of   $X_t$ do not see very fines things. 

An approximation with strongly stochastic errors  (cf above and \cite{bouleauRIMS}), particularly when the object to be approximated is a point of a functional space or a path of an erroneous stochastic process, may make visible only some properties of this path. \\

\noindent e) Up to now we have dealt with local Dirichlet forms admitting a carr\'e du champ operator, of course non local Dirichlet forms may possess also such an operator. It would mean to consider non local errors. This is not at all a crazy idea and such a concern appeared already in physics for example in Brillouin \cite{brillouin} Chap XV. About frequency rays of emission or absorption Brillouin distinguishes four possible cases of errors, one of which describes a very narrow ray with companion rays rather far from it with low probability. It seems that non local Dirichlet forms may be relevant in such cases for describing errors. 

Of course no differential calculus is available under such hypotheses for the propagation through calculations.\\

\noindent{\bf Concluding remarks.}
\\*With respect to the theory of Dirichlet forms, the error calculus \`a la Gauss looks quite like what the simply additive probability theory is with respect to the $\sigma$-additive theory. 

The error calculus by Dirichlet forms allows to perform computations on complex objects defined by limits -- as typically solutions to stochastic differential equations. But what yields it to engineers ? The question here is still very similar to the one we could ask about the probability calculus axiomatised by Kolmogorov in the framework of measure theory. This can bring much to engineers dealing with stochastic processes.

Not only finance in interested in processes\,! Especially all signal processings like Wiener and Kalman filtering, image improving, recognition, information transmission though channels \`a la Shannon, and particularly non linear treatments raise the question of an error calculus. We must get used to consider that any input process is accompanied with some accuracy defined by an error structure conditioning -- depending on the stochasticity and the intrinsic accuracy of the treatment -- the precision on the output expressed also by an error structure, so that a new treatment may be applied.


{\footnotesize
\begin{thebibliography}{20}

\bibitem{laplace}{\sc Laplace P. S.} "M\'emoire sur les probabilit\'es" {\it M\'emoirs de l'Acad\'emie royale des Sciences de Paris}, 1780
\OE uvres compl\`etes 9, pp. 227-332
\bibitem{legendre}{\sc Legendre A. M.} {\it Nouvelles m\'ethodes pour la d\'etermination des orbites des plan\`etes} Paris 1806.
\bibitem{adrain}{\sc Adrain R.} "Research concerning the probabilities of the errors which happens in making observations" Analyst math. museum, vol 1, n4, philadelphia 1808, 93-109.
\bibitem{laplace2}{\sc Laplace P. S.} {\it Th\'eorie analytique des probabilit\'es} Paris 1812.
\bibitem{gauss}{\sc Gauss C. F.} {\it Theoria motus corporum coelestium in sectionibus conicis solem ambientium }1809.
\bibitem{gauss0}{\sc Gauss C. F.} {\it Disquisitio de elementis ellipticis Palladis} 1811
\bibitem{laplace3}{\sc Laplace P. S.} {\it Essai philosophique sur les probabilit\'es} Paris, 1814.
\bibitem{cauchy}{\sc Cauchy A.} "M\'emoire [\ldots] pour que la plus grande de toutes les erreurs, abstraction faite du signe, devienne un maximum" 1\`ere classe de l'Institut 18 f\'ev.1814.
\bibitem{gauss}{\sc Gauss C. F.} "Bestimmung der Genauigkeit der Beobachtungen" {\it Lindenau und Bohnenberger Zeitschrift} 1816
\bibitem{gauss1}{\sc Gauss C. F.} {\it Theoria combinationis observationum erroribus minimis obnoxiae} 1821.
\bibitem{lacroix} {\sc Lacroix S. F.} {\it  Trait\'e \'el\'ementaire du calcul des probabilit\'es}, 1822
\bibitem{gauss2}{\sc Gauss C. F.} {\it Supplementum Theoriae Combinationis Observationum}   1828.
\bibitem{fourier}{\sc Fourier J.} "Second m\'emoire sur les r\'esultats moyens et sur les erreurs de mesure" Paris 1829.
\bibitem{bessel}{\sc Bessel F. W.} "Auszug aus einem Briefe des Herrn Prof. Bessel an den Herausgeber" {\it Astronomische Nachrichten}, vol 2, p.133, 1834
\bibitem{poisson}{\sc Poisson S. D.} {\it Recherches sur la probabilit\'e des jugements en mati\`ere criminelle et en mati\`ere civile, pr\'ec\'ed\'ees des r\`egles g\'en\'erales du calcul des probabilit\'es} 1837.
\bibitem{hagen}{\sc Hagen G.} {\it Gr\^undz\"uge der Wahrscheinlichkeitsrechnung} Berlin, 1837. 
\bibitem{bravi}{\sc Bravi G.} {\it Teorica e pratica del probabile}  Bergamo, 1840. 
\bibitem{cournot}{\sc Cournot A. A.} {\it Exposition de la th\'eorie des chances et des probabilit\'es} Paris, 1843. 
\bibitem{gerling}{\sc Gerling Ch. L.} {\it Die Angleichungsrechnungen der praktischen geometrie oder der Methode der kleinsten Quadrate} Hambourg, 1843. 
\bibitem{bravais}{\sc Bravais A.} {\it Analyse math\'ematique sur les probabilit\'es des erreurs de situation d'un point}. Paris 1844.
\bibitem{dirichlet}{\sc  Lejeune-Dirichlet G.} "Sur un moyen g\'en\'eral de v\'erifier l'expression du potentiel relatif \`a une masse quelconque, homog\`ene ou h\'et\'erog\`ene" {\it Crelle Journal f\"ur die reine und angewandte Mathematik}, Bd. 32, S. 80-94. (1846).
\bibitem{thomson}{\sc Thomson W.} "Note sur une \'equation aux diff\'erences partielles qui se pr\'esente dans plusieurs questions de Physique math\'ematique." {\it Journal de math\'ematiques pures et appliqu\'ees} 1re s\'erie, tome 12 (1847), p. 493-496.
\bibitem{riemann}{\sc Riemann B.} "Grundlagen f\"ur eine allgemeine Theorie der Funktionen einer ver\"anderlischen komplexen Gr\" osse (Inaugural Dissertation) 1851.
\bibitem{bienayme}{\sc Bienaym\'e I. J.} "Sur la probabilit\'e des erreurs" {\it Journal de math\'ematiques pures et appliqu\'ees} 1re s\'erie, tome 17 (1852), p. 33-78.
\bibitem{biver}{\sc Biver P. E.} "Th\'eorie analytique des moindres carr\'es" {\it Journal de math\'ematiques pures et appliqu\'ees }1re s\'erie, tome 18 (1853), p. 169-200.
\bibitem{cauchy1}{\sc Cauchy A.} "Sur la probabilit\'e des erreurs qui affectent des r\'esultats moyens d'observations de m\^eme nature" {\it C. R. Acad. Sci.}, t. XXXVII, p. 264 (1853)
\bibitem{cauchy2}{\sc Cauchy A.} "Sur la plus grande erreur \`a craindre dans un r\'esultat moyen, et sur le syst\`eme de facteurs qui rend cette plus grande erreur un minimum" {\it C. R. Acad. Sci.}, t. XXXVII, p. 326 (1853)
\bibitem{gauss3}{\sc Gauss C. F.} {\it M\'ethode des moindres carr\'es, m\'emoires sur la combinaison des observations} trad. J. Bertrand, Paris 1855.
\bibitem{riemann2}{\sc Riemann B.} ``Theorie der Abelschen Functionen" {\it Crelle Journal} 54, 115-155, 1857.
\bibitem{ameyer}{\sc Meyer A. } {\it Essai sur une exposition nouvelle de la th\'eorie analytique des probabilit\'es a posteriori.} Li\`ege, 1857.
\bibitem{airy}{\sc Airy G. B.} {\it On the algebraical and numerical theory of errors of observations and the combination of observations} Cambridge, 1861.
\bibitem{faadibruno}{\sc Faa di Bruno F.} {\it Trait\'e \'el\'ementaire du calcul des erreurs, avec des tables st\'er\'eotyp\'ees, ouvrage utile \`a ceux qui cultivent les sciences d'observation}  (1869).
\bibitem{liagre}{\sc Liagre J. B. J.} {\it Calcul des Porbabilit\'es et th\'eories des erreurs} Paris, Bruxelles 1879.
\bibitem{vonkries}{\sc von Kries, J.} {\it Die Prinzipien der Wahrscheinlichkeitsrechnung}, Freiburg 1886.
\bibitem{bertrand}{\sc Bertrand J.} {\it Calcul des probabilit\'es} 1888.
\vspace{.5cm}\\
{\large XXth century}

\bibitem{bachelier}{\sc Bachelier L.} "Th\'eorie des probabilit\'es continues", {\it J. de math\'ematiques pures et appliqu\'ees}, vol. 6, no 2, 1906, p. 259Ð327 
\bibitem{poincare}{\sc Poincar\'e, H.} {\it Calcul des Probabilit\'es} Gauthier-Villars, 1912.
\bibitem{palmer}{\sc A de Forest Palmer} {\it The Theory of Measurements} McGraw-Hill 1912
\bibitem{frechet}{\sc Fr\'echet M.} ``Remarque sur les probabilit\'es continues"{\it Bull. Sci. Math.} $2^e$ s\'erie, 45, (1921), 87-88.
\bibitem{borel} {\sc Borel E.} {\it Calcul des probabilit\'es} Paris, 1924.
\bibitem{fisher}{\sc Fisher R. A.}  {\it Statistical methods for research workers}. Genesis Publishing. (1925).
\bibitem{fisher2} {\sc Fisher R. A.} {\it Theory of Statistical Information}, Proc. Cambridge Philo. Soc. Vol XXII, Pt5,  1925.
\bibitem{hostinsky1}{\sc Hostinsk\'y B.} ``Sur la m\'ethode des fonctions arbitraires dans le calcul des probabilit\'es" {\it Acta Math.} 49, (1926), 95-113.
\bibitem{scarborough}{\sc Scarborough J. B.} "The invalidity of commonly used method for computing a certain probable error" {\it Proc. Nat. Acad. Sci} 15, 665 (1929).
\bibitem{schultz}{\sc Schultz H.} "Discussion" {\it J. of the American Statistical Association}, Vol. 24, No. 165, Supplement: Proceedings of the American Statistical Association (1929), pp. 86-89
\bibitem{kolmogorov0}{\sc Kolmogorov A. N.} "Ueber die analytischen Methoden in der Wahrscheinlichkeitsrecgnung", {\it Math. Annalen} 1931,
\bibitem{hostinsky2}{\sc Hostinsk\'y B.} {\it M\'ethodes g\'en\'erales de Calcul des Probabilit\'es}, Gauthier-Villars 1931.
\bibitem{birge}{\sc Birge R. T.} "The calculation of errors by the method of least squares"{\it Physical Review} V40, 1932.
\bibitem{hopf1}{\sc Hopf E.} ``On causality, statistics and probability" {\it J. of Math. and Physics} 18 (1934) 51-102.
\bibitem{hopf2} {\sc Hopf E.} ``\"{U}ber die Bedeutung der willk\"{u}rlichen Funktionen f\"{u}r die Wahrscheinlichkeitstheorie" {\it Jahresbericht der Deutschen Math. Vereinigung} XLVI, I, 9/12, 179-194, (1936).
\bibitem{hopf3}{\sc Hopf E.} ``Ein Verteilungsproblem bei dissipativen dynamischen Systemen" {\it Math. Ann.} 114, (1937), 161-186.
\bibitem{birge2}{\sc Birge R. T.} "The propagation of errors" {\it The Amer. Physics Teacher} V7, n6, 1939.
\bibitem{baret}{\sc Baret J.} {\it Calcul des erreurs} Paris 1939
\bibitem{hotelling}{\sc Hotelling H.} "Some new methods in matrix calculation" {\it Ann. Math. Statist.}, 14 (1943), pp. 1-34.
\bibitem{cartan}{\sc Cartan H.} ``Th\'eorie du potentiel newtonien'' {\it Bull. Soc. Math. de France} 73 (1945) 74-106.
\bibitem{cramer}{\sc Cram\'er,H.} {\it Mathematical methods of statistic} Princeton Univ. Press 1946.
\bibitem{vonneumann}{\sc Von Neumann J., Goldstine H. H. } "Numerical inverting of matrices of high order", {\it Bull. Amer. Math. Soc.}, 53 (1947), pp. 1021-1099.
\bibitem{turing}{\sc Turing A. M.} "Rounding-off errors in matrix processes " {\it Quart. J. Mech.}, (1948), pp.287-308.
\bibitem{vessiot}{\sc Vessiot E.,  Montel P.} {\it Cours de Math\'ematiques G\'en\'erales} Eyrolles 1947.
\bibitem{deny}{\sc Deny J.} "Les potentiels d'\'energie finie" {\it Acta Mathematica} 82, 107-183 (1950).
\bibitem{deny-bis}{\sc Deny J.} "Sur la d\'efinition de l'\'energie en th\'eorie du potentiel" {\it Ann. Inst. Fourier} 2, 83-99, 1950.
\bibitem{elkins1}{\sc Elkins Th.} "The second derivative method of gravity interpretation" {\it Society of Exploration Geophysicists} June 1950.
\bibitem{angot}{\sc Angot A.} {\it Compl\'ements de math\'ematiques \`a l'usage des ing\'enieurs de l'\'electrotechnique et des t\'el\'ecommunications} Paris 1952 
\bibitem{elkins2}{\sc Elkins Th.}"The effect of random errors in gravity data on second derivative values" {\it Geophysics}, January 1952, Vol. 17, No. 1 : pp. 70-88
\bibitem{parzen}{\sc Parzen E.} ``On estimation of a probability density function and mode" {\it Ann. Inst. Statist. Math. }6, 127-132, (1954)
\bibitem{ventsel}{\sc Wentzel A. D.} "Semi-groups of operators corresponding to a generalized differential operator of second order" {\it DAN SSSR}, Vol.111,2,1956,pp. 269-272.
\bibitem{hunt2}{\sc Hunt G. A.} "Markov processes and potentials II" {\it Illinois J. of Math.} 1957, n1, 316-369.
\bibitem{beurling-deny1}{\sc Beurling A., Deny J.} ``Espaces de Dirichlet I. Le cas \'el\'ementaire", {\it Acta Math.} 99, (1958) 203-224.
\bibitem{beurling-deny2}{\sc Beurling A., Deny J.}``Dirichlet spaces" {\it Proc. Nat. Acad. Sci. U.S.A.} 45 (1959) 208-215.
\bibitem{brillouin}{\sc Brillouin L.} {\it La Science et la Th\'eorie de l'Information} Paris, Masson 1959.
\bibitem{bartlett}{\sc Bartlett,M.S.} "Statistical estimation of density function" Sankhy\"a, sA, 25, 245-254 1963.
\bibitem{brelot1}{\sc Brelot M.} {\it Th\'eorie Classique du Potentiel} cours polycopi\'e CDU 3\`eme \'edition 1965.
\bibitem{murthy}{\sc Murthy,V. K.} "Estimation of probability density" {\it The Annals of Mathematical Statistics} ,6, 1027-103, 1965.
\bibitem{moore}{\sc Moore R. E.} {\it Interval Analysis} Prentice Hall 1966.
\bibitem{landau-lifchitz}{\sc Landau L. Lifchitz E.} {\it Physique Statistique} MIR 1967.
\bibitem{elkins3}{\sc Elkins Th.} "Cubical and Spherical Estimation of Multivariate Probability Density" {\it Journal of the American Statistical Association}, Vol. 63, No. 324 (Dec., 1968), pp. 1495- 1513.
\bibitem{deny1}{\sc Deny J.} ``M\'ethodes hilbertiennes en th\'eorie du potentiel'' Potential Theory, Centro Internazionale Matematico Estivo, Roma 1970.
\bibitem{wilkinson}{\sc Wilkinson J. H.} "Modern Error Analysis" {\it SIAM Review} Vol. 13, No. 4, October 1971.
 \bibitem{brelot}{\sc Brelot M.} ``Les \'etapes et les aspects multiples de la th\'eorie du potentiel'' {\it L'Enseignement Math\'ematique} XVIII, fasc.1, 1972.
 \bibitem{silverstein}{\sc Silverstein M. L. } {\it Symmetric Markov processes}, L. N. in Math. vol 426, Springer, 1974.
 \bibitem{klein}{\sc Klein H. A. } {\it The science of measurement, a historical survey}, Dover, 1974.
 \bibitem{deheuvels}{\sc Deheuvels P.} ``Estimation non param\'etrique de la densit\'e par histogrammes g\'en\'eralis\'es" {\it R. Statist. Appl.} Vol 25, f3, 1-24, (1977)
\bibitem{dellacherieAnglais}{\sc Dellacherie C.} and {\sc Meyer P.-A.} {\it Probability and Potential} North-Holland mathematics studies 29, 1978.
\bibitem{kolmogorov-yushkevich}{\sc Kolmogorov A. N., Yushkevich A. P.} {\it Mathematics of the 19th century (1978)}, Birkha\"user 1992.
\bibitem{pearson}{\sc Pearson K.} {\it The History of Statistics in the 17th \& 18th Centuries} Charles Griffin \& Co 1978.
\bibitem{lejan}{\sc Le Jan Y.} "Mesures associ\'ees \`a une forme de Dirichlet, applications" {\it Bull. Soc. Math. France} 106, 61-112, 1978.
\bibitem{fukushima0}{\sc Fukushima M.} "A decomposition of additive functionals of finite energy" {\it Nagoya Math. J.} 74, 1979, 137-168.
\bibitem{sheynin}{\sc Shetnin O. B.} "C. F. Gauss and the theory of errors" {\it Archive for History of Exact Sciences} 1979 Springer.
\bibitem{fukushima}{\sc Fukushima, M.}  {\it Dirichlet Forms and Symmetric Markov Processes} North-Holland 1980.
\bibitem{ikeda-watanabe}{\sc Ikeda N., Watanabe S.} {\it Stochastic Differential Equation and Diffusion Processes}, North-Holland, Koshanda 1981.
\bibitem{meyer81}{\sc Meyer P.-A.} "G\'eom\'etrie stochastique sans larmes" {\it S\'em. de Probabilit\'es XV} Lect. Notes in M. 850, Springer 1981.
\bibitem{temple}{\sc Temple G.} {\it 100 Years of Mathematics} Duckworth 1981.
\bibitem{bouleau1}{\sc Bouleau N.} "D\'ecomposition de l'\'energie par niveau de potentiel" {\it Colloque en l'honneur de Jacques Deny, Lect. Notes in M. 1096}, Springer(1984), http://hal.archives-ouvertes.fr/hal-00449195/fr/
\bibitem{ancona}{\sc Ancona A.} "L'\'energie et la th\'eorie du potentiel dans l'\oe uvre de Jacques Deny" {\it Colloque en l'honneur de Jacques Deny, Lect. Notes in M. 1096}, Springer(1984)
\bibitem{fuku-deny}{\sc Fukushima M.} "A Dirichlet form on the Wiener space and properties on Brownian motion ", {\it Colloque en l'honneur de Jacques Deny, Lect. Notes in M. 1096}, Springer(1984)
\bibitem{davydov-lifshits}{\sc Davydov Y. A., Lifshits M. A.} ``The stratification method in some probability problems'' {\it Prob. Theory, Math. Stat., Theor. Cybernetics}, 22,61-137, (1984).
\bibitem{lifshitz}{\sc Lifshitz M.A.} ``An application of the stratification method to the study of functionals of L\'evy processes and examples" {\it Theor. Probab. Appl.} \underline{29}, 753-764,(1984).
\bibitem{roeckner-wielens}{\sc R\"ockner M.} and {\sc Wielens N.} "Dirichlet forms --- Closability and change of speed measure", {\it Infinite Dimensional Analysis and Stochastic Processes}, Research Notes in M., Albeverio ed. Pitman 124, 119-144, (1985).
\bibitem{bouleau-hirsch1}{\sc Bouleau N.} and {\sc Hirsch F.}``Formes de Dirichlet g\'en\'erales et densit\'e des variables al\'eatoires r\'eelles sur l'espace de Wiener" {\it J. Funct. Analysis} 69, 2, 229-259,  (1986).
\bibitem{dellacherie4}{\sc Dellacherie C.} and {\sc Meyer P.-A.} {\it Probabilit\'es et Potentiel} Chap XII \`a XVI, Hermann 1987.
\bibitem{bichteler-gravereaux-jacod}{\sc Bichteler K., Gravereaux J.-B., Jacod J.} {\it Malliavin Calculus for Processes with Jumps} (1987).
\bibitem{bouleau1987}{\sc Bouleau N.} ``Autour de la variance comme forme de Dirichlet" in {\it S\'eminaire de Th\'eorie du Potentiel Paris No.8}, Lect. Notes in Math. 1235, Springer 1987.
\bibitem{wu}{\sc Wu L.} "Construction de l'op\'erateur de Malliavin sur l'espace de Poisson" {\it S\'em. Probabilit\'e XXI} Lect. Notes in M. 1247, Springer (1987).
\bibitem{jeudy}{\sc Jeudy L.} "Generalyzed variance-covariance propagation law formulae and application to explicit least-squares adjustments" {\it Bulletin g\'eod\'esique}
1988, Volume 62, Issue 2, pp 113-124
\bibitem{teunissen}{\sc Teunissen P. J. G.} "Nonlinearity and Least Squares" {\it CISM journ. ACSGC}, V42, n4, 1988
\bibitem{teunissen2}{\sc Teunissen P. J. G.} "First and Second Moments of non-Linear Least-Squares Estimators" {\it Bull. G\'eod.} 63 (1989) pp. 253-262.
\bibitem{coleman}{\sc Coleman H. W., Steele W. G.} {\it Experimentation and uncertainty analysis for engineers} Wiley 1989.
\bibitem{bouleau-lamberton}{\sc Bouleau N.} and {\sc Lamberton D.}
 "Residual risks and hedging strategies in Markovian markets", {\it Stochastic processes and their applications} 33,131-150, 1989.
 \bibitem{hangos}{\sc Hangos K. M., Leisztner L.} "Towards properly controlled analytical measurement methods" {\it J. of Automatic Chemistry} Vol. 11, No. 2 (March-April 1989), pp. 76-79
\bibitem{albeverio-roeckner}{\sc Albeverio S.} and {\sc R\"ockner M.} "Classical Dirichlet forms on topological vector spaces --- closability and a Cameron-Martin formula" {\it J. Funct. Analysis}, 88, 395-436,  (1990).
\bibitem{bouleau-hirsch2}{\sc Bouleau N.} and {\sc Hirsch F.} {\it Dirichlet Forms and Analysis on Wiener Space} De Gruyter (1991).
\bibitem{dellacherie5}{\sc Dellacherie C., Maisonneuve B.} and {\sc Meyer P.-A.} {\it Probabilit\'es et Potentiel} Chap XVII \`a XXIV, Hermann 1992.
\bibitem{engel}{\sc Engel E.} "A road to randomness in physical systems" Lect. notes in Stat. 71, Springer 1992.
\bibitem{ma-rockner1}{\sc Ma} and {\sc R\"ockner M.} {\it Introduction to the Theory of (Non-Symmetric) Dirichlet forms} Springer 1992.
\bibitem{albeverio-ma-roeckner}{\sc Albeverio S.,  Ma Z.-M., R\"ockner M.} ``Quasi-regular Dirichlet forms and Markov processes'', {\it J. Funct. Anal.} Ill (1993), 118-154.
\bibitem{coquio}{\sc Coquio A.} ``Formes de Dirichlet sur l'espace canonique de Poisson et application aux \'equations diff\'erentielles stochastiques" {\it Ann. Inst. Henri Poincar\'e} vol 19, n1, 1-36, (1993)
 \bibitem{bouleau-lepingle}{\sc Bouleau N.} and {\sc L\'epingle D.} {\it Numerical Methods for Stochastic Processes} Wiley (1994).
 \bibitem{fukushima-oshima-takeda}{\sc Fukushima M., Oshima Y.} and {\sc Takeda M.} {\it Dirichlet Forms and Symmetric Markov Processes} De Gruyter (1994).
 \bibitem{nualart}{\sc Nualart, D.} {\it Malliavin Calculus and Related Topics} Springer 1995.
\bibitem{bally-talay1}{\sc Bally V., Talay D.} ``The law of the Euler scheme for stochastic differential equations : I. Convergence rate of the distribution function", {\it Prob. Th. and Rel. Fields} vol 2 No2, 93-128 (1996).
\bibitem{bally-talay2}{\sc Bally V., Talay D.} ``The law of the Euler scheme for stochastic differential equations : II. Convergence rate of the density", {\it Monte Carlo Methods and Appl.} vol 104, No1, 43-80 (1996).
\bibitem{jacob}{\sc Jacob N.} {\it Pseudo-Differential Operators and Markov Processes}, Akademic Verlag (1996).
\bibitem{bogachev-et-al}{\sc Bogachev V., Roeckner M., Schmuland B.} ``Generalized Mehler semigroups and applications'' {\it Probability Theory and R.F.} V 105, N2, 193-225, (1996).
\bibitem{malliavin}{\sc  Malliavin P.}, {\it Sochastic Analysis}, Springer 1997.
 \bibitem{brelot2}{\sc Brelot M.} {\it Th\'eorie Classique du Potentiel} Ass. Laplace-Gauss Paris (1997)
 \bibitem{caballero-fernandez-nualart}{\sc  Caballero M. E., B. Fernandez, D. Nualart} ``Estimating densities and applications" {\it J. of Theoretical Probability} Vol 11, Nr3, (1998).
 \bibitem{silverman}{\sc Silverman B. W.} {\it Density Estimation for Statistics and Data Analysis} Chapman and Hall, 1998.
  \bibitem{villani}{\sc Villani C.} "Fisher Information Estimates for Boltzmann's Collision operator" {\it ?J. Math. Pures Appl.}, 77, 1998,p. 821-837.
\vspace{.5cm}\\
{\large XXIth century}

\bibitem{a-grothaus-k-r}{\sc Albeverio S., Grothaus M., Kondratiev Y. G., R\"ockner M.} ``Fluctuations in classical continuous systems" Univ. Bielefeld 2000.
\bibitem{ma-rockner2}{\sc Ma} and {\sc R\"ockner M.} ``Construction of diffusion on configuration spaces" {\it Osaka J. Math.} 37, 273-314,  (2000).
\bibitem{privault4}{\sc Privault N.} "Extended covariance identities and inequalities" {\it Statistics \& Prob. Letters}55, 247-255 (2001)
\bibitem{kohatsu-pettersson}{\sc Kohatsu-Higa A., Pettersson R.} ``Variance reduction methods for simulation of densities on Wiener space" SIAM J. Numer. Anal., 40(2), 431Ð450, 2002.
\bibitem{bouleau3}{\sc Bouleau N.} {\it Error Calculus for Finance and Physics, the Language of Dirichlet Forms}, De Gruyter (2003).
\bibitem{bouleau2003b}{\sc Bouleau N.} "Error calculus and path sensitivity in financial models"  {\it Mathematical Finance} vol 13 n$^0$ 1, 115-134, (Jan. 2003).
\bibitem{malliavin-thalmaier}{\sc Malliavin P., Thalmaier A.} ``Numerical error for SDE: Asymptotic expansion and hyperdistributions", {\it C. R. Acad. Sci. Paris}
ser. I 336 (2003) 851-856
 \bibitem{bouchard-ekeland-touzi}{\sc B. Bouchard, I. Ekeland, N. Touzi} ``On the Malliavin approach to Monte Carlo approximation of conditional expectations" {\it Finance Stochast.}
8, 45-71, (2004).
\bibitem{bouleauchorro}{\sc Bouleau N., Chorro Ch.}, "Error structures and parameter estimation" {\it C. R. Acad. Sci. Paris}, Ser. I 338, 305-310 (2004).
\bibitem{fitzsimmons}{\sc Fitzsimmons P.} "Superposition operators on Dirichlet spaces" Tohoku mathematical journal, 2004.
\bibitem{bouleauMC}{\sc Bouleau N.} ``Dirichlet forms in simulation"
 {\it Monte Carlo Methods and Appl.} Vol 11, n 4, 385-396, (2005).
 \bibitem{chorro}{\sc Chorro C.} ``Calculs d'erreur par formes de Dirichlet : liens avec l'information de Fisher et les th\'eor\`emes limites" Th\`ese, Univ. Paris I, 2005.
 \bibitem{bouleauD}{\sc Bouleau N.} "Th\'eor\`eme de Donsker et formes de Dirichlet"
{\it Bulletin des Sc. Math\'ematiques}, 129, (2005), 369-380
 \bibitem{chorro2}{\sc Chorro C.} ``Convergence in Dirichlet Law of Certain Stochastic
Integrals'' {\it Electr. J. of Probability} Vol 10, 1005-1025, (2005).
\bibitem{bouleauhal}{\sc Bouleau N.}``On some errors related to the graduation of measuring instruments" Oct. 2006
http://hal.archives-ouvertes.fr/hal-00105452/fr/
\bibitem{bouleau2006a}{\sc Bouleau N.}  "Bringing errors into focus" arXiv : 0705.0519
\bibitem{bouleau2006}{\sc Bouleau N.} ``When and how an error yields a Dirichlet form"
 {\it Journal of Functional Analysis}
Vol 240, Issue 2 , (2006)  445-494.
\bibitem{bally-et-al}{\sc Bally V., Bavouzet M.-P., Messaoud M.}``Intergration by parts formula for locally smooth laws and application to sensitivity computation" {\it The Annals of Appl. Prob.} Vol. 17, No. 1, 33Ð66, 2007.
 \bibitem{privault-wei}{\sc Privault N., Wei X.} ``Integration by parts for point processes and Monte Carlo simulation'' {\it J. of Appl. Prob.} 44, n3, 806-823, (2007).
  \bibitem{scotti2}{\sc Scotti S.}``Errors Theory using Dirichlet Forms, Linear Partial Differential Equations and Wavelets" arXiv:0708.1073, 2007
 {\bibitem{bouleau2006b}{\sc Bouleau N.} ``On error operators related to the arbitrary functions principle'' Jour. Functional Analysis 251, (2007) 325-345.}
 \bibitem{scotti}{\sc Scotti S.} {\it Applications de la Th\'eorie des Erreurs par Formes de Dirichlet}, Thesis Univ. Paris-Est, Scuola Normale Pisa, 2008. (http://pastel.paristech.org/4501/)
 \bibitem{scotti3}{\sc Scotti S.}``Perturbative Approach on Financial Markets" arXiv:0806.0287, 2008
\bibitem{bouleauRIMS}{\sc Bouleau N.} ``How to specify an approximate numerical result" {\it Proceedings of RIMS 2006 Workshop on Stochastic Analysis and Applications}, Kyoto 2007, {\it RIMS K\^oky\^uroku Bessatsu} B6, p39-53,  2008. http://hal.archives-ouvertes.fr/hal-00781414
\bibitem{scotti-lyvath}{\sc Scotti S., Ly Vath V.} ``Bid-ask Spread Modelling, a perturbation appoach"  finance-innovation.org (2008).
 \bibitem{bouleau-denis} {\sc Bouleau N.} and {\sc Denis L.} ``Energy image density property and the lent particle method for Poisson measures" {\it Jour. of Functional Analysis} 257 (2009) 1144-1174. available online: http//dx.doi.org/10.1016/j.jfa.2009.03.004
 \bibitem{privault5}{\sc Privault N.} {\it Stochastic Analysis in Discrete and Continuous Settings with Normal Martingales} Springer 2009.
  \bibitem{dembo-deuschel}{\sc Dembo A., Deuschel J.-D.}``Markovian perturbation, response and fluctuation dissipation theorem" arXiv:0710.4394v2 (2010).
  \bibitem{lim}{\sc Lim Th.-S.} "Quelques applications du contr\^ole stochastique aux risques de d\'efaut et de liquidit\'e" Thesis, Univ. Paris 7, 2010.
\bibitem{bouleau-denis2} {\sc Bouleau N.} and {\sc Denis L.} ``Application of the lent particle method to Poisson driven SDE's", {\it Probability Theory and Related Fields} Vol 151, Issue 3-4, pp 403-433, December 2011.
\bibitem{malicet-poly}{\sc Malicet D.} and {\sc Poly G.}``Propri\'et\'es de convergence dans les structures d'erreur" 2012 http://hal.archives-ouvertes.fr/hal-00608007
\bibitem{bouleaudenis}{\sc Bouleau N., Denis L.} {\it Dirichlet Forms Methods for Poisson Point Measures and L\'evy Processes,
with emphasis on creation-annihilation techniques} 265p, March 2013 to appear.


\end{thebibliography}

}
 \end{document}